\documentclass{amsart}
\usepackage[margin=1.5in]{geometry}

\usepackage[all]{xy}
\usepackage{tikz-cd}
\usepackage{enumerate}
\usepackage{amsfonts}
\usepackage{amssymb, amscd}
\usepackage{graphicx}
\usepackage{mathtools}
\usepackage{xcolor}
\usepackage{float}

\theoremstyle{definition}

\title{A short proof of Kneser's Theorem via transversals }
\author{Luis Montejano}

\begin{document}
\maketitle

\noindent \textbf{Kneser} \textbf{Theorem.}\cite{L}\ \textit{Whenever the }$k$%
\textit{-subsets of a }$n$\textit{-set are colored into }$\mathit{n-2k+1}$%
\textit{\ colors, then two disjoint }$k$-\textit{subsets end up having the
same color.}


For the proof we need the following definition.
A \textit{system of hyperplanes} $\chi $ in $\mathbb{R}^{d}$ consists of
choosing continuously one hyperplane of $\mathbb{R}^{d}$ in every direction.
That is, a system of hyperplanes in $\mathbb{R}^{d}$ is completely
determined by a function $\varphi :\mathbb{S}^{d-1}\rightarrow \mathbb{R},$
such that for every $x\in \mathbb{S}^{d-1},$ $\varphi (-x)=-\varphi (x).$
To see this, given $\varphi ,$ it is enough to choose, for every $x\in 
\mathbb{S}^{d-1},$ a hyperplane perpendicular to $x$ through $\varphi (x)x,$
which turns out to coincide with the hyperplane perpendicular to $-x$
through $-\varphi (-x)x,$ and vice versa.
The following lemma, that follows the spirit of Dolnikov, is crucial in the
proof of the Kneser Theorem.


\noindent \textbf{Lemma}. Given $d$ systems of hyperplanes $\chi
_{1},...,\chi _{d}$ in $\mathbb{R}^{d}$, there is a direction in which all
of them coincide.

\noindent \textbf{Proof}. For every $\mathbb{\chi }_{i},$ let $\varphi _{i}:%
\mathbb{S}^{d-1}\rightarrow \mathbb{R}$, such that $\varphi
_{i}(-x)=-\varphi _{i}(x)$, for every $x\in \mathbb{S}^{d-1},$ its
correspondent function. Let $\Psi :\mathbb{S}^{d-1}\rightarrow \mathbb{R}%
^{d-1},$ where $\Psi (x)=(\varphi _{1}(x)-\varphi _{d}(x),...,\varphi
_{d-1}(x)-\varphi _{d}(x))\in \mathbb{R}^{d-1},$ for every $x\in \mathbb{S}%
^{d-1}.$ Note that $\Psi (-x)=-\Psi (x).$ By the Borsuk Ulam Theorem, there
is $x_{0}\in \mathbb{S}^{d-1},$ such that $\Psi (x_{0})=0\in \mathbb{R}%
^{d-1}.$ That is, $\varphi _{i}(x_{0})=\varphi _{d}(x_{0}),$ for $%
i=1,...,d-1.$

\smallskip

\noindent \textbf{Proof of the Kneser Theorem}. Suppose not. Without loss of
generality, we may assume that our $n$-set is $V=\{x_{1},...x_{n}\}$ $%
\subset \mathbb{R}^{n-2k+1},$ a set of points in general position. A $k$%
-gone of $V$ is the convex closure of a $k$-subset of $V$. We assume that
the $k$-gons of $V$ are colored into $n-2k+1$ colors in such a way that two $%
k$-gons with the same color always intersect.

Let $L$ be a line through the origin. Choose a color, let say red. Then, the
projection of every red $k$-gone is a compact red interval contained in $L$.
Moreover, this collection of red intervals is mutually intersecting. Then,
by the Helly Theorem in the line, the intersection of all these red
intervals is again a compact interval. Take a hyperplane $H_{L},$
perpendicular to $L$, through the middle point of the intersection of all
red intervals. Clearly $H_{L}$ varies continuously with respect to $L,$
obtaining thus a red system of hyperplanes in $\mathbb{R}^{n-2k+1}.$
Similarly, we obtain a system of hyperplanes for every color. \ By our
Lemma, there is a direction in which all these $n-2k+1$ systems of
hyperplanes in $\mathbb{R}^{n-2k+1}$ coincide, but this implies that there
is a hyperplane $\Gamma $ transversal to all $k$-gons of $V.$

To the left of $\Gamma $ in $\mathbb{R}^{n-2k+1}$, we have less than $k$
points of $V$, otherwise $\Gamma $ would not be transversal to all $k$-gons
of $V.$ Similarly, to the right of $\Gamma $ in $\mathbb{R}^{n-2k+1},$ we
have less than $k$ points of $V$, so there are $n-2k+2$ or more points of $V$
in the hyperplane $\Gamma ,$ contradicting the fact that $V\subset \mathbb{R}%
^{n-2k+1}$ is a set of points in general position.

\end{document}